\documentclass{amsart}
\usepackage{amsfonts}
\usepackage{amssymb}
\usepackage{amsmath}
\usepackage{amsthm}

\newtheorem{theorem}{Theorem}

\newtheorem{lemma}{Lemma}

\theoremstyle{definition}
\newtheorem{definition}{Definition}

\begin{document}

\title{Asymptotically Conformal Fixed Points and Holomorphic Motions}

\author
{Yunping Jiang}
\address{Department of Mathematics\\
Queens College of the City University of New York\\
Flushing, NY 11367-1597\\
and\\
Department of Mathematics\\
Graduate School of the City University of New York\\
365 Fifth Avenue, New York, NY 10016}
\email[]{Yunping.Jiang@qc.cuny.edu}



\thanks{The research is partially supported by
NSF grants and PSC-CUNY awards and the Hundred Talents Program
from Academia Sinica.}

\begin{abstract}
The term integrable asymptotically conformal at a point for a
quasiconformal map defined on a domain is defined. Furthermore, we
prove that there is a normal form for this kind attracting or
repelling or super-attracting fixed point with the control
condition under a quasiconformal change of coordinate which is
also asymptotically conformal at this fixed point. The change of
coordinate is essentially unique. These results generalize
K\"onig's Theorem and B\"ottcher's Theorem in classical complex
analysis. The idea in proofs is new and uses holomorphic motion
theory and provides a new understanding of the inside mechanism of
these two famous theorems too.
\end{abstract}

\maketitle

\section{Introduction}

Two of the fundamental theorems in complex dynamical systems are
K\"onig's Theorem and B\"otthcher's Theorem in classical complex
analysis which were proved back to 1884~\cite{Kon} and
1904~\cite{Bot}, respectively, by using some well-known methods in
complex analysis. These theorems say that an attracting or
repelling or super-attracting fixed point of an analytic map can
be written into a normal form under a suitable conformal change of
coordinate. These theorems become two fundamental results in the
recent study of the dynamics of a polynomial or a rational map.

However, it becomes more and more clear in recent years that only
conformal changes of coordinate are not enough in the study of
many problems in dynamics and in geometry, for examples, in the
study of monotonicity of the entropy function for the family
$|x|^{3}+t$~\cite{MV}, in the study of deep points and
differentiability in hyperbolic $3$-manifolds~\cite[pp.32-34]{Mc}
(see also the end of \S3), and in the study of quasiconformal
structures on a $4$-manifold~\cite{DS}--in these studies,
quasiconformal changes of coordinate are appealed. The
quasiconformal changes of coordinate may still have asymptotical
conformality property just at one point but definitely not
conformal. (It is a big difference between asymptotically
conformal and conformal, see definition in Section 2.)

During the study of complex dynamical systems, a subject called
holomorphic motions becomes more and more interesting and useful.
The subject of holomorphic motions over the open unit disk shows
some interesting connections between classical complex analysis
and problems on moduli. This subject even becomes an interesting
branch in complex analysis~\cite{AM,BR,EM,Lie,Mit,ST}.

In this paper, we will use holomorphic motions over the open unit
disk to study the quasiconformal changes of coordinate which are
aymptotically conformal at one point. The paper is organized as
follows. In Section 2, we give an overview about holomorphic
motions and quasiconformal mapping theory. In Section 3, we define
an asymptotically conformal fixed point. Moreover, we define an
integrable asymptotically conformal fixed point. We then define an
attracting or repelling integrable asymptotically conformal fixed
point and the control condition. In Section 4, we prove one of our
main theorems in this paper:

\medskip
\begin{theorem}~\label{th1}
Let $f$ be a quasiconformal homeomorphism defined on a
neighborhood about $0$. Suppose $0$ is an attracting or repelling
integrable asymptotically conformal fixed point of $f$ with the
control condition. Then there is a quasiconformal homeomorphism
$\phi: \Delta_{\delta}\to \phi (\Delta_{\delta})\subset U$ from an
open disk of radius $\delta>0$ centered at $0$ into $U$ which is
asymptotically conformal at $0$ such that
$$
\phi^{-1} \circ f\circ \phi (z) =\lambda z, \quad z\in
\Delta_{\delta}.
$$
The conjugacy $\phi^{-1}$ is unique up to multiplication of a
constant.
\end{theorem}

\medskip
To present our idea clearly, we first use the same idea in the
proof of above theorem to give a new proof of K\"onig's Theorem in
classical complex analysis in Section 3. Then we prove
Theorem~\ref{th1} in the same section.

We define an asymptotically conformal super-attracting fixed point
in Section 3. In Section 5, we prove the other main theorem in
this paper:

\medskip
\begin{theorem}~\label{th2}
Let $g=f (z^{n})$ be a quasiregular map defined on a neighborhood
about $0$ for $n\geq 2$. Suppose $0$ is a super-attracting
integrable asymptotically conformal fixed point of $g$. Then there
is a quasiconformal homeomorphism $\phi: \Delta_{\delta}\to \phi
(\Delta_{\delta})\subset U$ from an open disk of radius $\delta>0$
centered at $0$ into $U$ which is asymptotically conformal at $0$
such that
$$
\phi^{-1} \circ g\circ \phi (z) =z^{n}, \quad z\in
\Delta_{\delta}.
$$
The conjugacy $\phi^{-1}$ is unique up to multiplication by
$(n-1)^{th}$-roots of the unit.
\end{theorem}

\medskip
Again, we will first give a new proof of B\"ottecher's Theorem in
classical complex analysis in Section 5. Then we prove
Theorem~\ref{th2} in the same section.

Our proofs in this paper use the ``holomorphic motion technique",
which we first used in~\cite{Ji}. Another place we used the
``holomorphic motion technique" is in the study of the Fatou
linearization and the quasiconformal rigidity for parabolic germs
in~\cite{Ji1}.

\vskip10pt \noindent {\bf Acknowledgment.} I would like to thank
Professor Weiyuan Qiu and my students Zhe Wang and Haifeng Chu to
help me to clarify several arguments and to fix mistakes and typos
in this paper. I would like also to thank Professor Sudeb Mitra
who explained to me several points in the development of the
measurable Riemann mapping theorem and holomorphic motions.

\section{Holomorphic Motions and Quasiconformal Maps}

In the study of complex analysis, the measurable Riemann mapping
theorem plays an important role. Consider the Riemann sphere
$\hat{\mathbb C}$. A measurable function $\mu$ on $\hat{\mathbb
C}$ is called a Beltrami coefficient if there is a constant $0\leq
k<1$ such that $\|\mu\|_{\infty} \leq k$, where
$\|\cdot\|_{\infty}$ means the $L^{\infty}$-norm of $\mu$ on
$\hat{\mathbb C}$. The equation
$$
H_{\overline{z}}=\mu H_{z}
$$
is called the Beltrami equation with the given Beltrami
coefficient $\mu$. The measurable Riemann mapping theorem says
that the Beltrami equation has a solution $H$ which is a
quasiconformal homeomorphism of $\hat{\mathbb C}$ whose
quasiconformal dilatation is less than or equal to
$K=(1+k)/(1-k)$. The study of the measurable Riemann mapping
theorem has a long history since Gauss considered in 1820's the
connection with the problem of finding isothermal coordinates for
a given surface. As early as 1938, Morrey~\cite{Mo} systematically
studied homeomorphic $L^{2}$-solutions of the Beltrami equation
(see~\cite{La1,La2}). But it took almost twenty years until in
1957 Bers~\cite{Be} observed that these solutions are
quasiconformal (refer to~\cite[pp. 24]{Le}). Finally the existence
of a solution to the Beltrami equation under the most general
possible circumstance, namely, for measurable $\mu$ with
$\|\mu\|_{\infty}<1$, was shown by Bojarski~\cite{Bo}. In this
generality the existence theorem is sometimes called the
measurable Riemann mapping theorem (refer to~\cite[pp. 10]{GL}.

If one only considers a normalized solution in the Beltrami
equation (a solution fixes $0$, $1$, and $\infty$), then $H$ is
unique, which is denoted as $H^{\mu}$. The solution $H^{\mu}$ is
expressed as a power series made up of compositions of singular
integral operators applied to the Beltrami equation on the Riemann
sphere. In this expression, if one considers $\mu$ as a variable,
then the solution $H^{\mu}$ depends on $\mu$ analytically. This
analytic dependence was emphasized by Ahlfors and Bers in their
1960 paper~\cite{AB} and is essential in determining a complex
structure for Teichm\"uller space (refer
to~\cite{Al,GL,Le,Li,Na}). Note that when $\mu\equiv 0$, $H^{0}$
is the identity map. A $1$-quasi-conformal map is conformal.
Twenty years later, due to the development of complex dynamics,
this analytic dependence presents an even more interesting
phenomenon called holomorphic motions as follows.

Let $\Delta_{r}=\{ c\in {\mathbb C} \;|\; |c|<r\}$ be the disk
centered at $0$ and of radius $r>0$. In particular, we use $\Delta$
to denote the unit disk. Given a Beltrami coefficient $\mu\not\equiv
0$, consider a family of Beltrami coefficients
$c\mu/\|\mu\|_{\infty}$ for $c\in \Delta$ and the family of
normalized solutions $H^{\frac{c\mu}{\|\mu\|_{\infty}}}$. Note that
$H^{\frac{c\mu}{\|\mu\|_{\infty}}}$ is a quasiconformal
homeomorphism whose quasiconformal dilatation is less than or equal
to $(1+|c|)/(1-|c|)$. Moreover, $H^{\frac{c\mu}{\|\mu\|_{\infty}}}$
is a family which is holomorphic on $c$. Consider a subset $E$ of
$\hat{\mathbb C}$ and its image
$E_{c}=H^{\frac{c\mu}{\|\mu\|_{\infty}}}(E)$. One can see that
$E_{c}$ moves holomorphically in $\hat{\mathbb C}$ when $c$ moves in
$\Delta$. That is, for any point $z\in E$,
$z(c)=H^{\frac{c\mu}{\|\mu\|_{\infty}}}(z)$ traces a holomorphic
path starting from $z$ as $c$ moves in the unit disk. Although $E$
may start out as smooth as a circle and although the points of $E$
move holomorphically, $E_{c}$ can be an interesting fractal with
fractional Hausdorff dimension for every $c\neq 0$ (see~\cite{GK}).

Surprisingly, the converse of the above fact is true too. This
starts from the famous $\lambda$-lemma of Ma\~n\'e, Sad, and
Sullivan~\cite{MSS} in complex dynamical systems. Let us start to
understand this fact by first defining holomorphic motions.

\medskip
\begin{definition}[Holomorphic Motions]~\label{hm} Let $E$ be a
subset of $\hat{\mathbb C}$. Let
$$
h (c, z) : \Delta_{r}\times E\to \hat{\mathbb C}
$$
be a map. Then $h$ is called a holomorphic motion of $E$
parametrized by $\Delta_{r}$ if
\begin{enumerate}
\item $h (0, z)=z$ for $z\in E$;
\item for any fixed $c\in
\Delta_{r}$, $h (c, \cdot): E\to \hat{\mathbb C}$ is injective;
\item
for any fixed $z$, $h (\cdot,z): \Delta_{r} \to \hat{\mathbb C}$
is holomorphic.
\end{enumerate}
\end{definition}

\medskip
For example, for a given Beltrami coefficient $\mu$,
$$
H(c, z)=H^{\frac{c\mu}{\|\mu\|_{\infty}}}(z): \Delta\times
\hat{\mathbb C}\to \hat{\mathbb C}
$$
is a holomorphic motion of $\hat{\mathbb C}$ parametrized by
$\Delta$.

Note that even continuity does not directly enter into the
definition; the only restriction is in the $c$ direction. However,
continuity is a consequence of the hypotheses from the proof of
the $\lambda$-lemma of Ma\~n\'e, Sad, and Sullivan~\cite[Theorem
2]{MSS}. Moreover, Ma\~n\'e, Sad, and Sullivan prove in~\cite{MSS}
that

\medskip
\begin{lemma}[$\lambda$-Lemma]~\label{ll}
A holomorphic motion $h(c,z): \Delta\times E\to \hat{\mathbb C}$
of a set $E\subset \hat{\mathbb C}$ parametrized by $\Delta$ can
be extended to a holomorphic motion $H(c,z): \Delta\times
\overline{E}\to \hat{\mathbb C}$ of the closure $\overline{E}$ of
$E$ parametrized by the same $\Delta$.
\end{lemma}

\medskip
Furthermore, Ma\~n\'e, Sad, and Sullivan show in~\cite{MSS} that
$H(c, \cdot): \overline{E}\to \hat{\mathbb C}$ satisfies the Pesin
property. In particular, when $\overline{E}$ is a closed domain,
this property can be described as the quasiconformal property. A
further study of this quasiconformal property is given by Sullivan
and Thurston~\cite{ST} and Bers and Royden~\cite{BR}.
In~\cite{ST}, Sullivan and Thurston prove that there is a
universal constant $a>0$ such that any holomorphic motion of any
set $E\subset \hat{\mathbb C}$ parametrized by the open unit disk
$\Delta$ can be extended to a holomorphic motion of $\hat{\mathbb
C}$ parametrized by $\Delta_{a}$. In~\cite{BR}, Bers and Royden
show, by using classical Teichm\"uller theory, that this constant
actually can be taken to be $1/3$. Moreover, in the same paper,
Bers and Royden show that in any holomorphic motion $H(c,z):
\Delta \times \hat{\mathbb C}\to \hat{\mathbb C}$, $H(c,\cdot):
\hat{\mathbb C}\to \hat{\mathbb C}$ is a quasiconformal
homeomorphism whose quasiconformal dilatation less than or equal
to $(1+|c|)/(1-|c|)$ for $c\in \Delta$. In the both
papers~\cite{ST,BR}, they expect $a=1$. This was eventually proved
by Slodkowski in~\cite{Sl}.

\medskip
\begin{theorem}[The Holomorphic Motion Theorem]~\label{st}
Suppose $h(c,z): \Delta\times E\to \hat{\mathbb C}$ is a
holomorphic motion of a closed subset $E$ of $\hat{\mathbb C}$
parameterized by the unit disk $\Delta$. Then there is a
holomorphic motion $H (c,z): \Delta \times \hat{\mathbb C}\to
\hat{\mathbb C}$ which extends $h(c,z): \Delta \times
E\to\hat{\mathbb C}$. Moreover, for any fixed $c\in \Delta$,
$H(c,\cdot): \hat{\mathbb C}\to\hat{\mathbb C}$ is a
quasiconformal homeomorphism whose quasiconformal dilatation
$$
K(H(c,\cdot)) \leq \frac{1+|c|}{1-|c|}.
$$
The Beltrami coefficient of $H(c, \cdot)$ given by
$$
\mu(c,z)=\frac{\partial H(c, z)}{\partial
\overline{z}}/\frac{\partial H(c, z)}{\partial z}
$$
is a holomorphic function from $\Delta$ into the unit ball of the
Banach space ${\mathcal L}^{\infty}({\mathbb C})$ of all
essentially bounded measurable functions on ${\mathbb C}$.
\end{theorem}

\medskip
Chirka gives a nice proof of Slodkowski's theorem by using some
results in functional analysis. The reader can find a complete
proof of the above holomorphic motion theorem in~\cite{GJW}
following the ideas in Chirka's proof~\cite{Ch} and in
Bers-Royden's proof~\cite{BR}. Moreover, some property of
infinitesimal holomorphic motions is discussed in~\cite{GJW}.

Holomorphic motions of a set $E\subset \hat{\mathbb C}$
parametrized by a connected complex manifold with a base point can
be also defined. They have many interesting relationships with the
Teichm\"uller space $T(E)$ of a closed set $E$ of the Riemann
sphere $\hat{\mathbb C}$ (refer to~\cite{Mit}).

In addition to the references we mentioned above, there is a
partial list of references~\cite{As,AM,Do,DHOP,EKK,EM,Lie,MS,Po}
about holomorphic motions and Teichm\"uller theory. The reader who
is interested in holomorphic motions may refer to those papers and
books.

\section{Integrable asymptotically conformal fixed points}

Let $f$ be a quasiconformal homeomorphism defined on a domain $U$
in the Riemann sphere $\hat{\mathbb C}$. Suppose $p$ is a point in
the $U$. Let $\Delta_{t} (p)$ denote the disk of radius $t>0$
centered at $p$. Let $\mu_{f}(z) = f_{\overline{z}}/f_{z}$ be the
complex dilatation of $f$ on $U$. Suppose $t_{0}>0$ be a number
such that $\Delta_{t_{0}}(p)\subset U$. Then for any $0<t\leq
t_{0}$, let $\omega_{f,p} (t) =\| \mu_{f}
|\Delta_{t}(p)\|_{\infty}$, where $\|\cdot\|_{\infty}$ means the
$L^{\infty}$ norm.

\medskip
\begin{definition}
We call $f$ asymptotically conformal at $p$ if
$$
\omega_{f,p} (t) \to 0 \quad \hbox{as} \quad t\to 0^+.
$$
Furthermore, we call $f$ integrable asymptotically conformal at
$p$ if
$$
\int_{0}^{t_{0}} \frac{\omega_{f,p} (s)}{s}ds <\infty.
$$
\end{definition}

\medskip
If $f$ is asymptotically conformal at $p$, then $f$ maps a tiny
circle centered at $p$ to an ellipse centered at $f(p)$ and,
moreover, the ratio of the long axis and the short axis tends to
$1$ as the radius of the tiny circle tends to $0$. But the map
still can fail to be differentiable at $p$ (refer to~\cite{GMRV}).
However, following Reshetnyak's 1978 paper~\cite[Theorem 1.1, pp.
204]{Re}), if $f$ is integrable asymptotically conformal at $p$,
then $f$ is differentiable and conformal at $p$, i.e., the limit
of $(f(z)-f(p))/(z-p)$ exists as $z$ goes to $p$. If, in addition,
$p$ is a fixed point of $f$, that is, $f(p)=p$, let
$$
\lambda= \lim_{z\to p} \frac{f(z)-f(p)}{z-p}
$$
and call it the multiplier of $f$ at $p$. We call $p$
\begin{itemize}
\item[i)] attracting if $0<|\lambda|<1$;
\item[ii)] repelling if $|\lambda|>1$;
\item[iii)] neutral if $|\lambda|=1$.
\end{itemize}
Correspondingly, we call $p$ an attracting, repelling, or neutral
integrable asymptotically conformal fixed point of $f$. By linear
changes of coordinate, we can assume that $p=f(p)=0$. We will keep
this assumption without loss of generality.

Let $g$ be a quasiregular map defined on a neighborhood $U$ of $0$
fixing $0$. Assume $g=f\circ q_{n}$ where $q_{n}(z)=z^{n}$, $n\geq
2$, and $f$ is a quasiconformal homeomorphism. We say $g$ is
integrable asymptotically conformal at $0$ if $f$ is integrable
asymptotically conformal at $0$ with nonzero multiplier
$$
\lambda =\lim_{z\to 0} \frac{f(z)}{z}.
$$
In this case, $0$ is called a super-attracting integrable
asymptotically conformal fixed point of $g$.

The following lemma will be useful in our proofs of
Theorems~\ref{th1} and~\ref{th2}.

\medskip
\begin{lemma}~\label{estimate}
Suppose $\omega (t)$ is an increasing function of $0<t\leq t_{0}$.
Suppose
$$
\int_{0}^{t_{0}} \frac{\omega (s)}{s} ds<\infty.
$$
Suppose $0<\sigma<1$ and $C>0$ are two constants. Let
$$
\tilde{\omega} (t) =\sum_{n=0}^{\infty} \omega (C\sigma^{n} t)
$$
for all $t>0$ such that $Ct\leq t_{0}$. Then
$$
\tilde{\omega}(t) \leq \omega(Ct) +\frac{1}{- \log \sigma}
\int_{0}^{Ct} \frac{\omega (s)}{s} ds.
$$
Moreover, $\tilde{\omega}(t) \to 0$ as $t\to 0^{+}$.
\end{lemma}

\medskip
\begin{proof} Since $\omega (t)$ is increasing for $t>0$, we have
$$
\tilde{\omega}(t) =\omega(Ct) +\sum_{n=1}^{\infty} \omega
(C\sigma^{n} t) \leq \omega(Ct) +\sum_{n=1}^{\infty}
\int_{n-1}^{n} \omega (C\sigma^{x} t) dx =\omega(Ct)
+\int_{0}^{\infty} \omega (C\sigma^{x} t) dx.
$$
Let $s=C\sigma^{x}t$. Then $ds =(\log \sigma) s dx$. We have that
$$
\int_{0}^{\infty} \omega (C\sigma^{x} t) dx = \frac{1}{-\log
\sigma}\int_{0}^{Ct} \frac{\omega (s)}{s} ds.
$$
\end{proof}

It is interesting to compare our integrable asympotically
conformal at a point to $C^{1+\alpha}$-conformal in McMullen's
book~\cite[pp.32]{Mc}. A homeomorphism $\phi(z)$ from a
neighborhood $U$ of ${\mathbb C}$ to another neighborhood $V$ of
${\mathbb C}$ is $C^{1+\alpha}$-conformal at $p\in U$ for some
$0<\alpha\leq 1$ if the complex derivative $\phi'(p)$ exists and
$$
\phi(z+p) = \phi(p) +\phi'(p) z + O(|z|^{1+\alpha})
$$
for all $z\in {\mathbb C}$ sufficiently small.

For a given Kleinian group $\Gamma$ preserving the upper-half
space ${\mathbb H}^{3}$ such that the $3$-manifold $M^{3}={\mathbb
H}^{3}/\Gamma$ has the bounded geometry, that is, its injectivity
radius is bounded above and below in its convex core. It is proven
that any quasiconformal conjugacy from $\Gamma$ to another
Kleinian group $\Gamma'$ is $C^{1+\alpha}$-conformal at every deep
point in the limit set $\Lambda$ of $\Gamma$ (see~\cite[Theorem
2.18]{Mc}). This theorem can be thought of as an extension of
Mostow rigidity at every deep point when the limit set is not the
whole sphere.

Suppose $\phi: U\to V$ is quasiconformal and
$C^{1+\alpha}$-conformal at $p$ for some $0<\alpha\leq 1$. By
linear changes of coordinate, we can assume $p=0$ and $\phi(p)=0$.
Then
$$
\phi (z) = \phi'(0) z + O(|z|^{1+\alpha})
$$
for $z\in {\mathbb C}$ sufficiently small. Since a quasiconformal
homeomorphism is differentiable almost everywhere, $F(z)=
O(|z|^{1+\alpha})$ is differentiable almost everywhere. Suppose
$F(z) = H(|z|^{1+\alpha})$. Then $H'(x)$ exists for almost every
$x$. Suppose $H'(x)$ is a bounded function. At every
differentiable point $z\not= 0$, we have that
$$
\Big| \frac{\partial \phi(z) }{\partial \overline{z}} \Big| =\Big|
\frac{\partial F(z)}{\partial \overline{z}} \Big| \leq C
|z^{\frac{\alpha+1}{2}} (\overline{z})^{\frac{\alpha-1}{2}}|
$$
and
$$
\Big|\frac{\partial \phi(z) }{\partial z}\Big| = \Big|\phi'(0) +
\frac{\partial F(z) }{\partial z}\Big| \geq |\phi'(0)| - C|
z^{\frac{\alpha-1}{2}} (\overline{z})^{\frac{\alpha+1}{2}}|
$$
where $C>0$ is a constant. Since $\phi'(0)\not=0$, we get that for
$z\in {\mathbb C}$ sufficiently small,
$$
|\mu_{\phi}(z)| =\Big|\frac{\partial \phi(z) }{\partial
\overline{z}}/\frac{\partial \phi(z) }{\partial z}\Big| \leq C'
|z|^{\alpha}
$$
where $C'>0$ is another constant. Then
$$
\omega_{\phi,0}(t)= \| \mu_{\phi}|\Delta_{t}\|_{\infty} \leq C'
t^{\alpha}.
$$
Thus $\phi$ is integrable asymptotically conformal at $0$.

If $\phi$ is quasiconformal in a neighborhood $U$ of $0$ fixing
$0$ and $C^{1+\alpha}$-conformal at $0$ for some $0<\alpha\leq 1$,
then it will automatically satisfy the control
condition~(\ref{control}) in the next section as follows.

Suppose $\lambda =\phi'(0)$ and suppose $0<|\lambda|<1$. (If
$|\lambda|>1$, then we consider $\phi^{-1}$.) Choose a constant
$0<a<1$ such that $a^{1+\alpha}<|\lambda|<a$. We can choose a
$\delta>0$ such that $\overline{\Delta}_{\delta}\subset U$ and
such that $|\phi (z)|\leq a|z|$ for any $z\in
\overline{\Delta}_{\delta}$. Then there is a constant $C_{0}>0$
such that, for any $|z|\leq \delta$,
$$
|\phi^{n+1} (z) -\lambda \phi^{n} (z) | \leq C_{0}
|\phi^{n}(z)|^{1+\alpha} \leq C_{0}|z|^{1+\alpha} a^{n(1+\alpha)}.
$$
Let $0<\tau = a^{1+\alpha}/|\lambda|<1$ and $C_{1}=
C_{0}\delta^{\alpha}/|\lambda|$. Then
$$
\Big| \frac{\phi^{n+1}(z)}{\lambda^{n+1}z}
-\frac{\phi^{n}(z)}{\lambda^{n}z} \Big| \leq C_{1} \tau^{n}.
$$
Thus $\{ h_{n}(z)= \phi^{n}(z)/(\lambda^{n}z)\}_{n=0}^{\infty}$ is
a uniform Cauchy sequence of continuous functions defined on
$\overline{\Delta}_{\delta}$. Furthermore, $h_{n}(0)=1$ for all
$n\geq 0$. Thus there is a constant $C>0$ such that
$$
C^{-1} \leq |h_{n}(z)| =\Big|\frac{\phi^{n}(z)}{\lambda^{n}z}\Big|
\leq C
$$
for all $z\in \overline{\Delta}_{\delta}$ and all $n\geq 0$ as
long as $\delta$ small enough. Therefore, if $\phi(z)$ is
quasiconformal and $C^{1+\alpha}$-conformal at $0$, then it
satisfies all assumptions in Theorem~\ref{th1}.

\medskip
\section{Linearization for integrable asymptotically conformal attracting or repelling fixed points}

One of the main results in this article is Theorem 1, which says
that if $f$ is a quasiconformal homeomorphism on a neighborhood
$U$ of $0$ and $0$ is an attracting or repelling integrable
asymptotically conformal fixed point with the multiplier
$\lambda$, $0<|\lambda |<1$ or $|\lambda|>1$ and with the control
condition, then $f$ can be written as a linear map $z \to \lambda
z$ under some quasiconformal change of coordinate which is also
asymptotically conformal at $0$. We only need to consider the
attracting case because that in the repelling case, we can
consider $f^{-1}$. In the attracting case, we say $f$ satisfies
the control condition if there are constants $\delta>0$ and $C>0$
such that
\begin{equation}~\label{control}
C^{-1} \leq \Big| \frac{f^{n}(z)}{\lambda^{n}z}\Big| \leq C
\end{equation}
for all $z\in \overline{\Delta}_{\delta}\subset U$ and all $n\geq
0$.

The result generalizes the famous K\"onig's Theorem in classical
analysis. Therefore, to present a clear idea about how we get
Theorem 1, we first use the same idea to give another proof of
K\"onig's Theorem, which is first given in~\cite{Ji}. The idea of
the new proof follows the viewpoint of holomorphic motions. For
the classical proof of K\"onig's Theorem, the reader may refer to
K\"onig's original paper~\cite{Kon} or most recent
books~\cite{CG,Mi}. Actually from the technical point of views,
our proof is more complicate and uses a sophistical result. But
from the conceptual point of views, our proof gives some inside
mechanism for the linearization of an attracting or a repelling
fixed point.

\medskip
\begin{theorem}[K\"onig's Theorem]~\label{kth}
Let $f(z)=\lambda z+\sum_{j=2}^{\infty} a_{j}z^{j}$ be an analytic
function defined on $\Delta_{r_{0}}$, $r_{0}>0$. Suppose
$0<|\lambda|<1$ or $|\lambda|>1$. Then there is a conformal map
$\phi: \Delta_{\delta} \to \phi(\Delta_{\delta})$ for some
$0<\delta <r_{0}$ such that
$$
\phi^{-1}\circ f\circ \phi (z) =\lambda z.
$$
The conjugacy $\phi^{-1}$ is unique up to multiplication of a
constant.
\end{theorem}

\medskip
\begin{proof}[A new proof of Theorem~\ref{kth}]
We only need to prove it for $0<|\lambda|<1$. In the case of
$|\lambda|>1$, we can consider $f^{-1}$.

First, we can find a $0<\delta<r_{0}$ such that
$$
|f(z)|< |z|, \quad z\in \overline{\Delta}_{\delta}
$$
and $f$ is injective on $\overline{\Delta}_{\delta}$. For every
$0<r\leq \delta$, let
$$
S_{r}=\{ z\in {\mathbb C}\;|\; |z|=r\}
$$ \
and
$$
T_{r}=|\lambda| S_{r}=\{ z\in {\mathbb C}\;|\; |z|=|\lambda| r\}.
$$
Denote $E=S_{r}\cup T_{r}$. Define
$$
\phi_{r} (z) =\left\{
\begin{array}{ll}
        z & z\in S_{r}\cr
        f(\frac{z}{\lambda}), & z\in T_{r}.
\end{array}\right.
$$
It is clear that
$$
\phi_{r}^{-1}\circ f\circ\phi_{r} (z)=\lambda z
$$
for $z\in S_{r}$.

Now write $\phi_{r}(z)=z\psi_{r} (z)$ for $z\in T_{r}$, where
$$
\psi (z) =1+\sum_{j=1}^{\infty} \frac{a_{j+1}}{\lambda^{j+1}}
z^{j}. $$ Define
$$
h_{r}(c, z) =\left\{
\begin{array}{ll}
z, & z\in S_{r}\cr z\psi (\frac{\delta cz}{r}), & z\in T_{r}
\end{array}\right.
: \Delta\times E\to \hat{\mathbb C}.
$$
Note that
$$
h_{r}(c,z) = z \psi \Big(
\frac{cz\delta}{r}\Big)=\frac{r}{c\delta} f\Big(
\frac{cz\delta}{r\lambda}\Big), \quad z\in T_{r}, \;\; c\neq 0.
$$
For each fixed $z\in E$, it is clear that $h (c, z)$ is a
holomorphic function of $c\in \Delta$. For each fixed
$c\in\Delta$, the restriction $h (c, \cdot)$ to $S_{r}$ and
$T_{r}$, respectively, are injective. Now we claim that their
images do not cross either. That is because for any $z\in T_{r}$,
$|z|=|\lambda|r$ and $|cz\delta|/|r\lambda| \leq \delta$, so
$$
|h(c, z)|=\Big| \frac{r}{c\delta}\Big| \Big| f\Big( \frac{
cz\delta}{r \lambda}\Big) \Big| <\Big| \frac{r}{c\delta}\Big|
\Big| \frac{cz\delta}{r\lambda}\Big| =r.
$$
Therefore, $h (c,z): \Delta\times E\to \hat{\mathbb C}$ is a
holomorphic motion because we also have $h(0, z) =z$ for all $z\in
E$. From Theorem 3, $h$ can be extended to a holomorphic motion
$H(c,z): \Delta\times \hat{\mathbb C}\to \hat{\mathbb C}$, and
moreover, for each fixed $c\in \Delta$, $H(c,\cdot):\hat{\mathbb
C}\to \hat{\mathbb C}$ is a quasiconformal homeomorphism whose
quasiconformal dilatation is less than or equal to
$(1+|c|)/(1-|c|)$. Now take $c_{r}=r/\delta$ and consider
$H(c_{r},\cdot)$. We have $H(c_{r}, \cdot)|E=\phi_{r}$. Let
$$
A_{r,j}= \{z\in {\mathbb C}\;|\; |\lambda|^{j+1}r\leq |z|\leq
|\lambda|^{j} r\}.
$$
We still use $\phi_{r}$ to denote $H(c_{r},\cdot)|A_{r, 0}$.

For an integer $k>0$, take $r=r_{k}=\delta |\lambda|^{k}$. Then
$$
\overline{\Delta}_{\delta} =\cup_{j=-k}^{\infty} A_{r,j}\cup
\{0\}.
$$
Extend $\phi_{r}$ to $\overline{\Delta}_{\delta}$, which we still
denote as $\phi_{r}$, as follows.
$$
\phi_{r}(z) = f^{j}(\phi_{r}(\lambda^{-j}z)), \quad z\in A_{r,j},
\quad j=-k, \cdots, -1, 0, 1, \cdots,
$$
and $\phi_{r}(0)=0$. Since $\phi_{r}|E$ is a conjugacy from $f$ to
$\lambda z$, $\phi_{r}$ is continuous on $\Delta_{\delta}$. Since
$f$ is conformal, $\phi_{r}$ is quasiconformal whose
quasiconformal dilatation is the same as that of $H(c_{r}, \cdot)$
on $A_{r,0}$. So the quasiconformal dilatation of $\phi_{r}$ on
$\Delta_{\delta}$ is less than or equal to $(1+c_r)/(1-c_r)$.
Furthermore,
$$
f (\phi_{r} (z)) =\phi_{r}(\lambda z), \quad z\in \Delta_{\delta}.
$$

Since $f^{-1}(z) =\lambda^{-1} z( 1+O(z))$, $f^{-k}(z)
=\lambda^{-k} z \prod_{i=0}^{k-1} (1+O(\lambda^{-i}z))$. Because
$|\lambda|^{-k}r_{k}=\delta$, the range of $\phi_{r_{k}}$ on
$\Delta_{\delta}$ is a Jordan domain bounded above from $\infty$
and below from $0$ uniformly on $k$. In addition, $0$ is fixed by
$\phi_{r_k}$ and the quasiconformal dilatations of the
$\phi_{r_k}$ are uniformly bounded. Therefore, the sequence $\{
\phi_{r_{k}}\}_{k=1}^{\infty}$ is in a compact set in the space of
all quasiconformal homeomorphisms on $\Delta_{\delta}$
(see~\cite{Al}). Let $\phi$ be a limiting map of this family. Then
we have
$$
f(\phi (z)) =\phi(\lambda z), \quad z\in \Delta_{\delta}.
$$
The quasiconformal dilatation of $\phi$ is less than or equal to
$(1+c_{r_{k}})/(1-c_{r_{k}})$ for all $k>0$. So $\phi$ is a
$1$-quasiconformal map, and thus is conformal. This is the proof
of the existence.

For the sake of completeness, we also provide the proof of
uniqueness but this is not new and the reader can find it
on~\cite{CG,Mi}. Suppose $\phi_{1}$ and $\phi_{2}$ are two
conjugacies such that
$$
\phi^{-1}_{1}\circ f\circ \phi_{1} (z) =\lambda z \quad
\hbox{and}\quad \phi^{-1}_{2}\circ f\circ \phi_{2} (z) =\lambda z,
\quad z\in \Delta_{\delta}.
$$
Then for $\Phi=\phi_{2}^{-1}\circ \phi_{1}$, we have $\Phi(\lambda
z)=\lambda \Phi(z)$. This implies that $\Phi'(\lambda z)=\Phi'(z)$
for any $z\in \Delta_{\delta}$. Thus
$\Phi'(z)=\Phi'(\lambda^{n}z)=\Phi(0)=0$. So
$\Phi(z)=\hbox{const}$ and $\phi_{2}^{-1}=\hbox{const.}\cdot
\phi_{1}^{-1}$.
\end{proof}

Now let us prove Theorem~\ref{th1}.

\begin{proof}[Proof of Theorem~\ref{th1}]
We need only to prove this theorem for attracting integrable
asymptotically conformal germs. In the case of repelling
integrable asymptotically conformal germs, we can consider
$f^{-1}$.

Let $\sigma=|\lambda|$. First, we can find a $\delta>0$ such that
$\overline{\Delta}_{\delta}\subset U$, $f$ is injective on
$\overline{\Delta}_{\delta}$,
$$
|f(z)|< |z|, \quad z\in \overline{\Delta}_{\delta},
$$
and the control condition~(\ref{control}) is held on
$\overline{\Delta}_{\delta}$.

For every $0<r\leq \delta$, let
$$
S_{r}=\{ z\in {\mathbb C}\;|\; |z|=r\}
$$ \
and
$$
T_{r}=\sigma S_{r}=\{ z\in {\mathbb C}\;|\; |z|=\sigma r\}.
$$
Denote $E=S_{r}\cup T_{r}$. Define
$$
\phi_{r} (z) =\left\{
\begin{array}{ll}
z & z\in S_{r}\cr f(\frac{z}{\lambda}), & z\in T_{r}
\end{array}
\right.
$$
It is clear that
$$
\phi_{r}^{-1}\circ f\circ\phi_{r} (z)=\lambda z
$$
for $z\in S_{r}$.

Now write $\phi(z) =f(\frac{z}{\lambda})$ defined on $\Delta_{r}$.
Suppose $\phi (r)=\tau_{r}$. Extend $\phi$ to $\hat{\mathbb C}$ by
quasiconformal reflection with respect to $S_{r}$ and
$\phi(S_{r})$ (see~\cite{Al}). We still denote this extended map
as $\phi$. Let $\nu=\phi_{\overline{z}}/\phi_{z}$ be the complex
dilatation of the extended $\phi$. Then
$$
a(r)=\| \nu\|_{\infty} = O(\|\mu
|\Delta_{\sigma^{-1}r}\|_{\infty})=O(\omega (\sigma^{-1}r)).
$$
Consider $\nu_{c}= ca_{0}a(r)^{-1} \nu$ and the unique solution
$\phi_{c} =\phi^{\nu_{c}}$ that maps $0$, $r$, and $\infty$ to $0$,
$\tau_{r}$, and $\infty$, respectively. Here $a_{0}$ is a constant
independent of $r$ such that $|\phi_{c}(z)| <r$ for all $|z|\leq
\sigma r$ and $|c|<1$. (Since $\phi_{c}$ can be written as a power
series in $c$ and $\|\nu_{c}\|\to 0$ uniformly as $r\to 0$ , such an
$a_{0}$ exists.) Then $\phi_{c}$ holomorphically depends on $c\in
\Delta$. Define
$$
\phi_{r} (c, z) =\left\{
\begin{array}{ll}
z & z\in S_{r}, \cr \phi_{c} (z), & z\in T_{r}.
\end{array}\right.
$$
It is a holomorphic motion from $\Delta\times E \to \hat{\mathbb
C}$. From Theorem~\ref{st}, $\phi_{r} (c, z)$ can be extended to a
holomorphic motion from $\Delta\times \hat{\mathbb C} \to
\hat{\mathbb C}$, which we still denote by $\phi_{r} (c, z)$, such
that the quasiconformal dilatation of $\phi_{r} (c, \cdot)$ is
less than or equal to $(1+|c|)/(1-|c|)$. In particular when
$c_{r}=a_{0}^{-1}a(r)$, $\phi_{r} (c_{r},z)|E =\phi_{r}$. Let
$$
A_{r,j}= \{z\in {\mathbb C}\;|\; \sigma^{j+1}r\leq |z|\leq
\sigma^{j} r\}.
$$
We still use $\phi_{r}$ to denote $\phi_{r}(c_{r}, \cdot)|A_{r,
0}$. For an integer $k>0$, take $r=r_{k}=\delta \sigma^{k}$. Then
$$
\Delta_{\delta} =\cup_{j=-k}^{\infty} A_{r,j}\cup \{0\}.
$$
Extend $\phi_{r}$ to $\Delta_{\delta}$, which we still denote as
$\phi_{r}$, by
$$
\phi_{r}(z) = f^{j}(\phi_{r}(\lambda^{-j}z)), \quad z\in A_{r,j},
\quad j=-k, \cdots, -1, 0, 1, \cdots,
$$
and $\phi_{r}(0)=0$. Since $\phi_{r}|E$ is a conjugacy from $f$ to
$\lambda z$, $\phi_{r}$ is continuous on $\Delta_{\delta}$.

Next we need to estimate the quasiconformal constant of $\phi_{r}$
on $\Delta_{\delta}$. We will use the following formula (refer
to~\cite{Al}): If $F$ and $G$ are two quasiconformal maps with the
complex dilatations $\mu_{F}$ and $\mu_{G}$. Then the composition
map $G\circ F$ has the complex dilatation
\begin{equation}~\label{composition}
\mu_{G\circ F} = \frac{\mu_{F}+\gamma \mu_{G}\circ F}{1+
\overline{\mu}_{F} \gamma \mu_{G}\circ F}, \quad \hbox{where}\quad
\gamma = \frac{\overline{F}_{\overline{z}}}{F_{z}}.
\end{equation}
Thus
$$
\|\mu_{G\circ F}\|_{\infty} \leq (\| \mu_{F}\|_{\infty} +
\|\mu_{G}\circ F\|_{\infty}) (1- \| \mu_{F}\|_{\infty} \|
\mu_{G}\circ F\|_{\infty})^{-1}.
$$

Let $\omega (t) =\omega_{f,0}(t)$. Suppose $C>0$ is the constant
in the control condition~(\ref{control}). Suppose, in the
beginning of the proof, we pick $\delta$ small such that $\omega
(C\delta) <1$. From Lemma~\ref{estimate},
$$
K_{0}=\tilde{\omega} (\delta) =\sum_{n=0}^{\infty} \omega
(C\sigma^{n} \delta) <\infty
$$
is a convergent series. Thus the product
$$
K_{1} =\prod_{n=0}^{\infty} \Big(1-\omega (C\sigma^{n}
\delta)\Big)^{-1}<\infty
$$
is also convergent.

Let $\mu (z) =\mu_{\phi_{r}}(z)$ for $z\in \Delta_{\delta}$.
Remember that $r=\sigma^{k}\delta$. For $z\in A_{r,0}$, $|\mu(z)|
\leq c_{r}$. For $z\in A_{r, -j}$, $1 \leq j\leq k$, $\phi_{r}(z)
= f^{-j} (\phi_{r}(\lambda^{j}z))$. Let $g_{i} (z) =
f^{-j+i}(\phi_{r}(\lambda^{j}z))$ for $0\leq i\leq j$. Let
$w=\lambda^{j}z$. Then $\sigma r\leq |w|\leq r$. This implies that
$|\phi_{r}(w)| \leq r$ and
$$
|g_{i}(z)| \leq C\sigma^{-j+i} r = C\sigma^{k-j+i} \delta
$$
for all $0\leq i\leq j$.  Note that
$$
|\mu_{f^{-1}}| =|\mu_{f}\circ f^{-1}|.
$$
By the compsoition formula~(\ref{composition}),
$$
|\mu(z)|=|\mu_{g_{0}}(z)| \leq \Big(|\mu_{g_{1}}(z)|+|\mu_{f^{-1}}
(g_{1}(z))|\Big) \Big(1-|\mu_{g_{1}}(z)||\mu_{f^{-1}}
(g_{1}(z))|\Big)^{-1}
$$
$$
= \Big(|\mu_{g_{1}}(z)|+|\mu_{f} (g_{0}(z))|\Big)
\Big(1-|\mu_{g_{1}}(z)||\mu_{f} (g_{0}(z))|\Big)^{-1}
$$
$$
\leq |\mu_{g_{1}}(z)|\Big(1-\omega (C\sigma^{k-j}\delta)\Big)^{-1}
+\omega (C\sigma^{k-j}\delta) \Big(1-\omega
(C\sigma^{k-j}\delta)\Big)^{-1}.
$$
Inductively, we get
$$
|\mu_{g_{i}}(z)| \leq |\mu_{g_{i+1}}(z)|\Big(1-\omega
(C\sigma^{k-j+i}\delta)\Big)^{-1} +\omega (C\sigma^{k-j+i}\delta)
\Big(1-\omega (C\sigma^{k-j+i}\delta)\Big)^{-1}
$$
for $0\leq i\leq j$. So
$$
|\mu(z)| \leq c_{r} \prod_{l=1}^{j}
\Big(1-\omega(C\sigma^{k-l}\delta)\Big)^{-1} +\sum_{i=1}^{j}\omega
(C\sigma^{k-i}\delta) \prod_{l=k-j}^{k-i}
\Big(1-\omega(C\sigma^{l}\delta)\Big)^{-1}
$$
$$
\leq K_{1}c_{r} +K_{1} \sum_{i=1}^{j}\omega (C\sigma^{k-i}\delta)
\leq K_{1} (c_{r} +\tilde{\omega}(\delta)) \leq K_{1} (1+K_{0}).
$$

For $z\in A_{r, j}$, $1 \leq j<\infty$, $\phi_{r}(z) = f^{j}
(\phi_{r}(\lambda^{-j}z))$. Let $h_{i} (z) =
f^{j-i}(\phi_{r}(\lambda^{-j}z))$ for $0\leq i\leq j$. Let
$w=\lambda^{-j}z$. Then $\sigma r\leq |w|\leq r$. This implies
that $|\phi_{r}(w)| \leq r$ and
$$
|h_{i}(z)| \leq C\sigma^{j-i} r = C\sigma^{k+j-i} \delta
$$
for all $0\leq i\leq j$. By the compsoition
formula~(\ref{composition}),
$$
|\mu(z)|=|\mu_{h_{0}}(z)| \leq \Big(|\mu_{h_{1}}(z)|+|\mu_{f}
(h_{1}(z))|\Big) \Big(1-|\mu_{h_{1}}||\mu_{f}
(h_{1}(z))|\Big)^{-1}
$$
$$
\leq |\mu_{h_{1}}(z)|\Big(1-\omega (C\sigma^{k+j-1}
\delta)\Big)^{-1} +\omega (C\sigma^{k+j-1} \delta) \Big(1-\omega
(C\sigma^{k+j-1} \delta)\Big)^{-1}.
$$
Inductively, we get
$$
|\mu_{h_{i}}(z)| \leq |\mu_{h_{i+1}}(z)|\Big(1-\omega
(C\sigma^{k+j-i-1}\delta)\Big)^{-1} +\omega
(C\sigma^{k+j-i-1}\delta) \Big(1-\omega
(C\sigma^{k+j-i-1}\delta)\Big)^{-1}.
$$
So
$$
|\mu(z)| \leq c_{r} \prod_{l=1}^{j}
\Big(1-\omega(C\sigma^{k+j-l}\delta)\Big)^{-1}
+\sum_{i=1}^{j}\omega (C\sigma^{k+j-i}\delta) \prod_{l=1}^{i}
\Big(1-\omega(C\sigma^{k+j-l}\delta)\Big)^{-1}
$$
$$
\leq K_{1} c_{r} + K_{1} \sum_{i=1}^{j}\omega
(C\sigma^{k+j-i}\delta) \leq K_{1} (c_{r}
+\tilde{\omega}(\sigma^{k}\delta)) =K_{1} (c_{r}
+\tilde{\omega}(r))\leq K_{1} (1+K_{0}).
$$

Let $k=K_{1}(1+K_{0})$ and $K=(1+k)/(1-k)$. Then $\{
\phi_{r_{k}}\}_{k=1}^{\infty}$ is uniformly $K$-quasiconformal.
Consider $B_{r}=\Delta_{\delta}\setminus
\Delta_{r}=\cup_{j=-k}^{-1}A_{r,j}$ and $\phi_{r}
(B_{r})=\cup_{j=-k}^{-1}\phi_{r}(A_{r,j})$ for any $r=r_{k}$. Both
of the annulli have the same inner circle $S_{r}$. Thus the ratio
of the modulus of $\phi_{r}(B_{r})$ and the modulus of $B_{r}$ is
controlled by two constants from below and above (independent of
$r$ but only depends on $K$). Therefore, the range of $\phi_{r}$
on $\Delta_{\delta}$ is a Jordan domain bounded above from
$\infty$ and below from $0$ uniformly in $0<r=r_{k}\leq \delta$.
Since, additionally, $0$ is fixed by any element in this sequence,
the family $\{ \phi_{r_{k}}\}_{k=1}^{\infty}$ is in a compact set
in the space of all $K$-quasiconformal homeomorphisms on
$\Delta_{\delta}$ (see~\cite{Al}). Let $\phi$ be a limit mapping
of this family. Then we have
$$
f(\phi (z)) =\phi(\lambda z), \quad z\in \Delta_{\delta}.
$$

Similar to the arguments above, the complex dilatation of
$\phi_{r}(z)$ on disk $\Delta_{\tilde{r}}$ is controlled by
$K_{1}(c_{\tilde{r}}+ \tilde{\omega}(\tilde{r}))$ for any
$r=r_{k}\leq \tilde{r}$. So the complex dilatation of $\phi$ on
$\Delta_{\tilde{r}}$ is also controlled by $K_{1}(c_{\tilde{r}}+
\tilde{\omega}(\tilde{r}))\to 0$ as $\tilde{r}\to 0$. Thus
$\phi(z)$ is asymptotically conformal at $0$ and the proof of
existence is completed.

Suppose $\phi_{1}$ and $\phi_{2}$ are two asymptotically conformal
conjugacies such that
$$
\phi^{-1}_{1}\circ f\circ \phi_{1} (z) =\lambda z \quad
\hbox{and}\quad \phi^{-1}_{2}\circ f\circ \phi_{2} (z) =\lambda z,
\quad z\in \Delta_{\delta}.
$$
Then for $\Phi=\phi_{2}^{-1}\circ \phi_{1}$, we have $\Phi(\lambda
z)=\lambda \Phi(z)$. This implies that the complex dilatation
$\mu_{\Phi} (z) = \mu_{\Phi} (\lambda z)$, a.e.. This in turn
implies that $\mu=0$ a.e. in $\Delta_{\delta}$ and thus $\Phi$ is
conformal. Furthermore, $\Phi(z) =az$ for some $a\not=0$. This is
the uniqueness.
\end{proof}

\medskip
\section{Normal forms for integrable asymptotically conformal
super-attracting fixed points}

The other main result in this article is Theorem~\ref{th2}, which
says that if $g=f (z^{n})$ is a quasiregular map and $0$ is an
integrable asymptotically conformal super-attracting fixed point,
then $g$ can be written into the normal form $z: \to z^n$ under
some quasiconformal change of coordinate which is asymptotically
conformal at $0$. The result generalizes the famous B\"ottcher's
Theorem in classical analysis. Again, to present a clear idea
about how we get Theorem~\ref{th2}, we first use the same idea to
give another proof of B\"ottcher's Theorem, which is first given
in~\cite{Ji}. The idea of the new proof follows the viewpoint of
holomorphic motions. For the classical proof of B\"ottcher's
Theorem, the reader may refer to B\"ottcher's original
paper~\cite{Bot} or most recent books~\cite{CG,Mi}. Actually from
the technical point of views, our proof is more complicate and
uses a sophistical result. But from the conceptual point of views,
our proof gives some inside mechanism of the normal form for a
super-attracting fixed point. The idea of the proof is basically
the same as that in the previous section, but the actual proof is
little bit different. The reason is that in the previous case, $f$
is a homeomorphism so we can iterate both forward and backward,
but in Theorem~\ref{th2} or B\"ottcher's Theorem, $g$ is not a
homeomorphism.

\medskip
\begin{theorem}[B\"ottcher's Theorem]~\label{bth} Suppose $g(z)
=\sum_{j=n}^{\infty}a_{j}z^{j}$, $a_{n}\neq 0$, $n\geq 2$, is
analytic on a disk $\Delta_{\delta_{0}}$, $\delta_{0}>0$. Then
there exists a conformal map $\phi: \Delta_{\delta}\to
\phi(\Delta_{\delta})$ for some $\delta>0$ such that
$$
\phi^{-1}\circ g\circ \phi (z) =z^{n}, \quad z\in \Delta_{\delta}.
$$
The conjugacy $\phi^{-1}$ is unique up to multiplication by
$(n-1)^{th}$-roots of the unit.
\end{theorem}

\medskip
\begin{proof}[A new proof of B\"ottcher's Theorem]
Conjugating by $z\to bz$, we can assume $a_{n}=1$, i.e.,
$$
g(z) =z^{n} +\sum_{j=n+1}^{\infty} a_{j}z^{j}.
$$

We use $\Delta_{r}^{*}=\Delta_{r}\setminus \{0\}$ to mean a
punctured disk of radius $r>0$. Write
$$
g(z)=z^{n}(1+\sum_{j=1}^{\infty} a_{j+n}z^{j}).
$$
Assume $0<\delta_{1}<\min\{ 1/2, \delta_{0}/2\}$ is small enough
such that
$$
1+\sum_{j=1}^{\infty} a_{j+n}z^{j} \neq 0 \quad \hbox{and} \quad
\frac{1}{\hbox{$\root n \of  {|1+\sum_{j=1}^{\infty}
a_{j+n}z^{j}}|$}} \geq \frac{1}{2}, \quad z\in
\Delta_{2\delta_{1}}.
$$
Then $g: \Delta_{2\delta_{1}}^{*}\to g(\Delta_{2\delta_{1}}^{*})$
is a covering map of degree $n$.

Let $0<\delta < \delta_{1}$ be a fixed number such that
$g^{-1}(\Delta_{\delta}) \subset \Delta_{\delta_{1}}$. Since
$$
z\to z^{n}: \Delta_{\hbox{$\root n \of {\delta}$}}^{*} \to
\Delta_{\delta}^{*} \quad \hbox{and}\quad  g:
g^{-1}(\Delta_{\delta}^{*}) \to \Delta_{\delta}^{*}
$$
are both of covering maps of degree $n$, the identity map of
$\Delta_{\delta}$ can be lifted to a holomorphic diffeomorphism
$$
h: \Delta_{\hbox{$\root n \of {\delta}$}}^{*}\to
g^{-1}(\Delta_{\delta}^{*}),
$$
i.e., $h$ is a map such that the diagram
$$
\begin{array}{ccc}
\Delta_{\hbox{$\root n \of {\delta}$}}^{*} & {\buildrel h \over
\longrightarrow}& g^{-1}(\Delta_{\delta}^{*})\cr \downarrow z\to
z^{n} &         &\downarrow g\cr \Delta_{\delta}^{*} &{\buildrel
\hbox{id} \over \longrightarrow}&
   \Delta_{\delta}^{*}
\end{array}
$$
commutes. We pick the lift so that
$$
h (z) =z \Big( 1+\sum_{j=2}^{\infty} b_{j}z^{j-1}\Big) =z\psi(z).
$$
From
$$
g(h(z))=z^{n}, \quad z\in \Delta_{\hbox{$\root n \of
{\delta}$}}^{*},
$$
we get
$$
|h(z)| =\frac{|z|}{\hbox{$\root n \of
{|1+\sum_{j=1}a_{n+j}(h(z))^{j}|}$}} \geq \frac{|z|}{2}.
$$

For any
$$
0<r\leq \min \Big\{ \Big( \frac{1}{2}\Big)^{\frac{n}{(n-1)}},
\delta^{n}\Big\},
$$
let $S_{r}=\{ z\in {\mathbb C}\; |\; |z|=r\}$ and $T_{r}=\{ z\in
{\mathbb C}\;|\; |z|=\hbox{$\root n \of r$}\}$. Consider the set
$E=S_{r}\cup T_{r}$ and the map
$$
\phi_{r}(z) =\left\{
\begin{array}{ll}
z, & z\in S_{r}\cr
z\psi(z), & z\in T_{r}.
\end{array}\right.
$$
Define
$$
h_{r}(c, z) =\left\{
\begin{array}{ll}
z, & z\in S_{r}\cr z\psi \Big( \frac{cz}{\hbox{$\root n \of
r$}}\Big), & z\in T_{r}
\end{array}\right.
:\Delta\times E\to \hat{\mathbb C}.
$$
Note that
$$
z\psi\Big( \frac{cz}{\hbox{$\root n \of r$}}\Big)
=\frac{\hbox{$\root n \of r$}}{c} h \Big( \frac{cz}{\hbox{$\root n
\of r$}}\Big), \quad z\in T_{r},\; c\not= 0.
$$
This implies that
$$
|h_{r}(c,z)| =\frac{\hbox{$\root n \of r$}}{|c|} h\Big(
\frac{cz}{\hbox{$\root n \of r$}}\Big) \Big| \geq
\frac{\hbox{$\root n \of r$}}{|c|} \frac{|cz|}{2 \hbox{$\root n
\of r$}} \geq \frac{\hbox{$\root n \of r$}}{2}>r, \quad z\in
T_{r}.
$$
So images of $S_{r}$ and $T_{r}$ under $h_{r}(c,z)$ do not cross
each other.

Now let us check $h_{r}(c,z)$ is a holomorphic motion. First
$h_{r}(0,z)=z$ for $z\in E$. For fixed $x\in E$, $h_{r}(c, z)$ is
holomorphic on $c\in \Delta$. For fixed $c\in \Delta$,
$h_{r}(c,z)$ restricted to $S_{r}$ and $T_{r}$, respectively, are
injective. But the images of $S_{r}$ and $T_{r}$ under
$h_{r}(c,z)$ do not cross each other. So $h_{r}(c, z)$ is
injective on $E$.  Thus
$$
h_{r}(c,z): \Delta\times E\to \hat{\mathbb C}
$$
is a holomorphic motion. By Theorem~\ref{st}, it can be extended
to a holomorphic motion
$$
H_{r}(c, z): \Delta\times \hat{\mathbb C}\to \hat{\mathbb C}.
$$
And moreover, for each $c\in \Delta$, $H_{r}(c,\cdot)$ is a
quasiconformal map whose quasiconformal dilatation satisfies
$$
K(H_{r}(c, \cdot)) \leq \frac{1+|c|}{1-|c|}.
$$

Now consider $H_{r}(\hbox{$\root n \of r$}, \cdot)$. It is a
quasiconformal map with quasiconformal constant
$$
K_{r}\leq \frac{1+\hbox{$\root n \of r$}}{1-\hbox{$\root n \of
r$}}.
$$
Let
$$
A_{r,j}=\{z\in {\mathbb C}\;|\; \hbox{$\root n^{j} \of r$} \leq
|z|\leq \hbox{$\root n^{j+1} \of r$} \},\quad j=0, 1, 2, \cdots.
$$
Consider the restriction $\phi_{r,0}=H_{r}(\hbox{$\root n \of r$},
\cdot)|A_{r, 0}$. It is an extension of $\phi_{r}$, i.e.,
$\phi_{r,0}|E=\phi_{r}$.

Let $\tilde{A}_{r,0}$ be the annulus bounded by $S_{r}$ and
$g^{-1}(S_{r})$ and define $\tilde{A}_{r, j} =g^{-j}
(\tilde{A}_{r, 0})$, $j\geq 0$. Since $z\to z^{n}: A_{r,1}\to
A_{r, 0}$ and $g: \tilde{A}_{r,1}\to \tilde{A}_{r,0}$ are both
covering maps of degree $n$, so $\phi_{r,0}$ can be lifted to a
quasiconformal map $\phi_{r, 1}: A_{r, 1}\to \tilde{A}_{r, 1}$,
i.e., the following diagram
$$
\begin{array}{ccc}
A_{r,1} &{\buildrel \phi_{r,1} \over \longrightarrow}&
\tilde{A}_{r,1}\cr \downarrow z\to z^{n} &         &\downarrow
g\cr A_{r,0} &{\buildrel \phi_{r,0} \over \longrightarrow}&
\tilde{A}_{r,0}
\end{array}
$$
commutes. We pick the lift $\phi_{r, 1}$ such that it agrees with
$\phi_{r,0}$ on $T_{r}$. The quasiconformal dilatation of
$\phi_{r,1}$ is less than or equal to $K_{r}$.

For an integer $k>0$, take $r=r_{k}=\delta^{n^k}$. Inductively, we
can define a sequence of $K_{r}$-quasiconformal maps $\{\phi_{r,
j}\}_{j=0}^{k}$ such that
$$
\begin{array}{ccc}
A_{r,j} &{\buildrel \phi_{r,j} \over \longrightarrow}&
\tilde{A}_{r,j}\cr \downarrow z\to z^{n} &         &\downarrow
g\cr A_{r,j-1} &{\buildrel \phi_{r,j-1} \over \longrightarrow}&
\tilde{A}_{r,j-1}
\end{array}
$$
commutes and $\phi_{r,j}$ and $\phi_{r, j-1}$ agree on the common
boundary of $A_{r,j}$ and $A_{r, j-1}$. Note that
$$
\overline{\Delta}_{\delta}=\Delta_{r}\cup \cup_{j=0}^{k} A_{r,j}.
$$
Now we can define a quasiconformal map, which we still denote by
$\phi_{r}$ as follows.
$$
\phi_{r} (z) =\left\{
\begin{array}{lll} z, & z\in \Delta_{r}; &\cr
\phi_{r, j}, & z\in A_{r, j},& j=0, 1, \cdots, k.
\end{array}\right.
$$
The quasiconformal dilatation of $\phi_{r}$ on $\Delta_{\delta}$
is less than or equal to $K_{r}$ and
$$
g(\phi_{r}(z)) =\phi_{r}(z^{n}), \quad z\in \cup_{j=1}^{k}
A_{r,j}.
$$

Since $g(z) =z^{n}( 1+O(z))$, $g^{k}(z) = z^{n^{k}}
\prod_{i=0}^{k-1} (1+O(z^{n^i}))$. Because ${\root n^{k} \of
r_{k}}=\delta$, the range of $\phi_{r_{k}}$ on $\Delta_{\delta}$
is a Jordan domain bounded above from $\infty$ and below from $0$
uniformly in $k$. In addition, $0$ is fixed by $\phi_{r_k}$ and
the quasiconformal dilatations of the $\phi_{r_k}$ are uniformly
bounded in $k$. Therefore, the sequence $\{
\phi_{r_{k}}\}_{k=1}^{\infty}$ is in a compact set in the space of
all quasiconformal homeomorphisms on $\Delta_{\delta}$
(see~\cite{Al}). Let $\phi$ be a limiting map of this family. Then
we have
$$
g(\phi (z)) =\phi(z^{n}), \quad z\in \Delta_{\delta}.
$$
Since the quasiconformal dilatation of $\phi$ is less than or
equal to $(1+{\root n \of r_{k}}) /(1-{\root n \of r_{k}})$ for
all $k>0$, it follows that $\phi$ is a $1$-quasiconformal map, and
thus conformal. This is the proof of the existence.

Suppose $\phi_{1}$ and $\phi_{2}$ are two conjugacies such that
$$
\phi^{-1}_{1}\circ g\circ \phi_{1} (z) =z^{n} \quad
\hbox{and}\quad \phi^{-1}_{2}\circ g\circ \phi_{2} (z) =z^{n},
\quad z\in \Delta_{\delta}.
$$
For
$$
\Phi (z) =\phi_{2}^{-1}\circ \phi_{1} (z) =\sum_{j=1}^{\infty}
a_{j}z^{j},
$$
we have $\Phi( z^{n})=(\Phi(z))^{n}$. This implies
$a_{1}^{n}=a_{1}$ and $a_{j}=0$ for $j\geq 2$. Since $a_{1}\neq
0$, we have $a_{1}^{n-1}=1$ and
$\phi_{2}^{-1}=a_{1}\phi_{1}^{-1}$. This is the uniqueness.
\end{proof}

\medskip

We now prove Theorem~\ref{th2}. The proof follows almost the same
footsteps of those of Theorem~\ref{th1} and Theorem~\ref{bth}.

\medskip
\begin{proof}[Proof of Theorem~\ref{th2}.]
Let $g=f\circ q_{n}$, $n\geq 2$. Conjugating by $z\to bz$, we can
assume $f'(0)=\lim_{|z|\to 0} f(z)/z=1$.

We use $\Delta_{r}^{*}=\Delta_{r}\setminus \{0\}$ to mean a
punctured disk of radius $r>0$. There is a $0<\delta_{1}<1$ such
that $g: \Delta_{2\delta_{1}}^{*}\to g(\Delta_{2\delta_{1}}^{*})$
is a covering map of degree $n$.

Let $0<\delta < \delta_{1}$ be a fixed number such that
$g^{-1}(\Delta_{\delta}) \subset \Delta_{\delta_{1}}$. Since
$$
z\to z^{n}: \Delta_{\hbox{$\root n \of {\delta}$}}^{*} \to
\Delta_{\delta}^{*} \quad \hbox{and}\quad  g:
g^{-1}(\Delta_{\delta}^{*}) \to \Delta_{\delta}^{*}
$$
are both of covering maps of degree $n$, the identity map of
$\Delta_{\delta}$ can be lifted to a homeomorphism
$$
h: \Delta_{\hbox{$\root n \of {\delta}$}}^{*}\to
g^{-1}(\Delta_{\delta}^{*}).
$$
Furthermore, $h$ is a quasiconformal map and integrable
asymptotically conformal at $0$ such that the diagram
$$
\begin{array}{ccc}
\Delta_{\hbox{$\root n \of {\delta}$}}^{*} & {\buildrel h \over
\longrightarrow}& g^{-1}(\Delta_{\delta}^{*})\cr \downarrow z\to
z^{n} &         &\downarrow g\cr \Delta_{\delta}^{*} &{\buildrel
\hbox{id} \over \longrightarrow}&
   \Delta_{\delta}^{*}
\end{array}
$$
commutes. We pick the lift so that
$$
h'(0) =\lim_{z\to 0} \frac{h(z)}{z}=1.
$$
These can be seen from the equation
$$
g(h(z))=z^{n}, \quad z\in \Delta_{\hbox{$\root n \of
{\delta}$}}^{*}.
$$

For any $0<r\leq \delta$, let $S_{r}=\{ z\in {\mathbb C}\; |\;
|z|=r\}$ and $T_{r}=\{ z\in {\mathbb C}\;|\; |z|=\hbox{$\root n
\of r$}\}$. Consider the set $E=S_{r}\cup T_{r}$ and the map
$$
\phi_{r}(z) =\left\{
\begin{array}{ll}
z, & z\in S_{r}\cr h(z), & z\in T_{r}.
\end{array}\right.
$$

It is clear that
$$
 g (\phi_{r} (z))=\phi_{r} (z^n)
$$
for $z\in T_{r}$.

Extend $h$ to $\hat{\mathbb C}$ by quasiconformal reflection with
respect to $S_{r}$ and $\phi(S_{r})$ (see~\cite{Al}). We still
denote this extended map as $\phi$. Let
$\nu=\phi_{\overline{z}}/\phi_{z}$ be the complex dilatation of
the extended $\phi$. Then
$$
a(r)= \|\nu\|_{\infty} =O(\| \mu |\Delta_{{\root n \of
r}}\|_{\infty})=O(\omega ({\root n \of r})).
$$
Assume $h({\root n \of r}) =\tau_{r}$. Consider $\nu_{c}=
ca_{0}a(r)^{-1} \nu$ and the unique solution $\phi_{c}
=\phi^{\nu_{c}}$ that maps $0$, $r$, and $\infty$ to $0$,
$\tau_{r}$, and $\infty$, respectively. Here $a_{0}$ is a constant
independent of $r$ such that $|\phi_{c}(z)| >r$ for all $|z|=
\sqrt[n]{r}$ and $|c|<1$. (Since $\phi_{c}$ can be written as a
power series in $\nu_{c}$ and $\|\nu_{c}\|\to 0$ uniformly as
$r\to 0$ , such an $a_{0}$ exists.) Then $\phi_{c}$
holomorphically depends on $c\in \Delta$. Define
$$
\phi_{r} (c, z) =\left\{
\begin{array}{ll}
z & z\in S_{r}\cr \phi_{c} (z), & z\in T_{r}.
\end{array}\right.
$$
It is a holomorphic motion from $\Delta\times E \to \hat{\mathbb
C}$. From Theorem~\ref{st}, $\phi (c, z)$ can be extended to a
holomorphic motion from $\Delta\times \hat{\mathbb C} \to
\hat{\mathbb C}$, which we still denote by $\phi (c, z)$, such
that the quasiconformal dilatation of $\phi (c, \cdot)$ is less
than or equal to $(1+|c|)/(1-|c|)$. In particular when
$c_{r}=a_{0}^{-1}a(r)$, $\phi_{r} (c_{r},z)|E =\phi_{r}$. We still
use $\phi_{r}$ to denote $\phi_{r}(c_{r}, \cdot)|A_{r, 0}$.

For an integer $k>0$, take $r=r_{k}= \delta^{n^k}$. Let
$$
A_{r,j}= \{z\in {\mathbb C}\;|\; {\root n^{j} \of r}\leq |z|\leq
{\root n^{j+1} \of r}\}
$$
for $0\leq j\leq k-1$. Then
$$
\Delta_{\delta} =\Delta_{r}\cup \cup_{j=0}^{k-1} A_{r,j}.
$$
Let $\phi_{r} (z) =z$ for $z\in \Delta_{r}$ and extend $\phi_{r}$
to $\cup_{j=0}^{k-1} A_{r,j}$ by lifting. Then we get a
homeomorphism on $\Delta_{\delta}$, which we still denote as
$\phi_{r}$ (refer to the proof of Theorem~\ref{th1}). Formally we
can use the following formula to define $\phi_{r}$,
$$
\phi_{r}(z) = g^{-j}(\phi_{r}(z^{n^{j}})), \quad z\in A_{r,j},
\quad j=0, \cdots, k-1,
$$
and $\phi_{r}(z)=z$ for $z\in \Delta_{r}$. Since $\phi_{r}|E$ is a
conjugacy from $g$ to $q_{n}(z) =z^{n}$, $\phi_{r}$ is continuous
on $\Delta_{\delta}$.

Let $\omega (t) =\omega_{f,0}(t)$ for $0<t\leq \delta$. Suppose
$C>0$ is a constant such that
$$
C^{-1} \leq \Big|\frac{f(z)}{z}\Big| \leq C
$$
for $z\in \overline{\Delta}_{\delta}$. Suppose, in the beginning
of the proof, we pick $\delta$ small such that $\omega (C\delta)
<1$. Let $0<\delta <\sigma<1$ be a fixed constant. From
Lemma~\ref{estimate},
$$
K_{0}=\tilde{\omega} (\delta) =\sum_{n=0}^{\infty} \omega
(C\sigma^{n} \delta) <\infty
$$
is a convergent series. Thus the product
$$
K_{1} =\prod_{n=0}^{\infty} \Big(1-\omega (C\sigma^{n}
\delta)\Big)^{-1}<\infty
$$
is also convergent.

Using the similar argument to that in the proof of
Theorem~\ref{th1}, we can obtain that the complex dilatation
$\mu(z)=\mu_{\phi_{r}}(z)$ over $\Delta_{\delta}$ can be
controlled by
$$
|\mu(z)| \leq K_{1} (c_{r} +\tilde{\omega}(\delta))\leq K_{1}
(1+K_{0})
$$
for $z\in \Delta_{\delta}$ and
$$
|\mu(z)| \leq K_{1} (c_{\tilde{r}} +\tilde{\omega}(\tilde{r}))
$$
for $z\in \Delta_{\tilde{r}}$ and all $r=r_{k}\leq \tilde{r}$.

Let $k=K_{1}(1+K_{0})$ and $K=(1+k)/(1-k)$. Then $\{
\phi_{r_{k}}\}_{k=1}^{\infty}$ is uniformly $K$-quasiconformal.
Consider $B_{r}=\Delta_{\delta}\setminus
\Delta_{r}=\cup_{j=0}^{k-1}A_{r,j}$ and $\phi_{r}
(B_{r})=\cup_{j=0}^{k-1}\phi_{r}(A_{r,j})$ for any $r=r_{k}$. Both
of the annulli have the same inner circle $S_{r}$. Thus the ratio
of the modulus of $\phi_{r}(B_{r})$ and the modulus of $B_{r}$ is
controlled by two constants from below and above (independent of
$r$ but only depends on $K$). Therefore, the range of $\phi_{r}$
on $\Delta_{\delta}$ is a Jordan domain bounded above from
$\infty$ and below from $0$ uniformly in $0<r=r_{k}\leq \delta$.
Since, additionally, $0$ is fixed by any element in this sequence,
the sequence $\{ \phi_{r_{k}}\}_{k=1}^{\infty}$ is in a compact
set in the space of all $K$-quasiconformal homeomorphisms on
$\Delta_{\delta}$ (see~\cite{Al}). Let $\phi$ be a limit mapping
of this family. Then we have
$$
g(\phi (z)) =\phi(z^{n}), \quad z\in \Delta_{\delta}.
$$

Similar to the arguments in Theorem~\ref{th1}, the complex
dilatation of $\phi_{r}(z)$ on disk $\Delta_{\tilde{r}}$ is
controlled by $K_{1}(c_{\tilde{r}}+ \tilde{\omega}(\tilde{r}))$
for any $r=r_{k}\leq \tilde{r}$. So the complex dilatation of
$\phi$ on $\Delta_{\tilde{r}}$ is also controlled by
$K_{1}(c_{\tilde{r}}+ \tilde{\omega}(\tilde{r}))\to 0$ as
$\tilde{r}\to 0$. Thus $\phi(z)$ is asymptotically conformal at
$0$. The proof of existence is completed.

Suppose $\phi_{1}$ and $\phi_{2}$ are two asymptotically conformal
conjugacies such that
$$
\phi^{-1}_{1}\circ g\circ \phi_{1} (z) =z^{n} \quad
\hbox{and}\quad \phi^{-1}_{2}\circ g\circ \phi_{2} (z) =z^{n},
\quad z\in \Delta_{\delta}.
$$
Then for $\Phi=\phi_{2}^{-1}\circ \phi_{1}$, we have
$\Phi(z^{n})=(\Phi(z))^{n}$. This implies that the complex
dilatation $\|\mu_{\Phi} (z)\| = \|\mu_{\Phi} (z^{n})\|$, a. e..
This in turn implies that $\mu=0$ a.e. in $\Delta_{\delta}$ and
thus $\Phi$ is conformal, and therefore, $\Phi (z) =az$ with
$a^{n-1}=1$. This is the uniqueness.
\end{proof}

\bibliographystyle{amsalpha}

\end{document}